\documentclass[12pt]{article}

\usepackage[T2A]{fontenc}
\usepackage[utf8]{inputenc}

\usepackage{amsfonts}
\usepackage{amssymb}
\usepackage{amsmath}
\usepackage{amsthm}

\def \le {\leqslant}
\def \ge {\geqslant}

\theoremstyle{plain}

\begin{document}
\begin{Huge}
 \centerline{A note on two linear forms}
\end{Huge}
\vskip+0.5cm
\centerline{by {\bf Klaus Moshchevitin}\footnote{ Research is supported by RFBR grant No.12-01-00681-a and by the grant of Russian Government, project 11.
G34.31.0053.
  }  }
\vskip+1.5cm

{\bf 1. Diophantine exponents.}

Let $\theta_1,\theta_2$ be  real numbers  such that
\begin{equation}\label{ine}
1,\theta_1,\theta_2\,\,\,
\text{are linearly independent over }\,\,\mathbb{Z}.
\end{equation}
We
consider linear form
$$
L({\bf x})= x_0+x_1\theta_1+x_2\theta_2,
\,\,\,\,
{\bf x} = (x_0,x_1,x_2) \in \mathbb{Z}^3.
$$
By $|{\bf z}|$ we denote the Euclidean length of a vector $ {\bf z} =(z_0,z_1,z_2) \in \mathbb{R}^3$.
Let
 \begin{equation}\label{2}
\hat{\omega} = \hat{\omega}(\theta_1,\theta_2) =
\sup \left\{\gamma:\,\,
\limsup_{t\to\infty}\,\,\left(
t^\gamma \min_{0<|{\bf x}|\le t}{|L({\bf x})|}\right) <\infty \right\}
\end{equation}
be the uniform Diophantine exponent
for the linear form $L$.

We consider  another linear form
$P({\bf x})$. The main result of the present paper is as follows.

{\bf Theorem 1.}\,\,{\it
Suppose that linear forms
$L(x)$ and $P(x)$ are independent and the exponent $\hat{\omega}$ for the form $L$ are defined in (\ref{2}).
Then for the Diophantine exponent
$$
\omega_{LP} =
\sup \left\{\gamma:\,\,
\text{there exist infinitely many}\,\,
{\bf x}\in \mathbb{Z}^3\,\,\text{such that}\,\,
|L({\bf x})|\le |P({\bf x})| \cdot |{\bf x}|^{-\gamma}
 \right\}
$$
we have a lower bound
$$
\omega_{LP} \ge
 \hat{\omega}^2 - \hat{\omega} +1.
$$
}

{\bf Remark.}
\,\, Of course in the definition (\ref{2}) and in  Theorem 1 instead of the Euclidean norm $|{\bf x}|$ we may consider the value $ \max_{i=1,2} |x_j|$ as it was done  by the most of authors.

Consider a  real $\theta$ which is not a rational number and not a quadratic irrationality.
Define
$$
\omega_* =\omega_*(\theta)=\sup\{
\gamma:\,\,\text{there exist infinitely many algebraic numbers}\,\, \xi\,\,\text{ of degree}\,\, \le 2\,\,
$$
$$
\text{such that}\,\,
|\theta - \xi | \le H(\xi)^{-\gamma}\}
$$
(here $H(\xi )$ is the maximal value of the absolute values of the coefficients for canonical polynomial to $\xi$).
Then for linear forms
$$
L({\bf x})= x_0+x_1\theta+x_2\theta^2,\,\,\,\,  P({\bf x})= x_1+2x_2\theta
$$
one has
\begin{equation}\label{diff}
\omega_* \ge \omega_{LP}.
\end{equation}

So Theorem 1 immediately leads to the following corollary.

{\bf Theorem 2.}\,\,{\it For a real
$\theta$ which is not a rational number and not a quadratic irrationality
one has
\begin{equation}\label{t2}
\omega_* \ge
 \hat{\omega}^2 - \hat{\omega} +1
\end{equation}
with  $\hat{\omega} = \hat{\omega} (\theta, \theta^2).$
}

{\bf 2. Some history.}

In 1967 H. Davenport and W. Schmidt \cite{DS1}  (see also Ch. 8 from Schmidt's book \cite{SB}) proved
that for any two independent linear forms $L,P$  there exist infinitely many integer points ${\bf x}$ such that
$$
|L({\bf x})| \le C |P({\bf x})|\, |{\bf x}|^{-3},
$$
with  a positive constant $C$ depending on the coefficients of forms $L,P$. From this result they deduced that for
any real $\theta$ which is not a rational number and not a quadratic irrationality the inequality
$$
|\theta - \xi |\le C_1 H(\xi)^{-3}
$$
has infinitely many solutions in algebraic $\xi$ of degree $\le 2$.

We see that  for any two pairs of forms one has $\omega_{LP} \ge 3$. But form the Minkowski convex body theorem
it follows that  under the condition (\ref{ine}) one has $\hat{\omega} \ge 2$. Moreover
$$
\min_{\hat{\omega}\ge 2}   (\hat{\omega}^2 - \hat{\omega} +1)  = 3.
$$
So our Theorems 1,2 may be considered as generalizations of Davenport-Schmidt's results.

Later Davenport and Schmidt
generalized their theorems to the case of several linear forms \cite{DS2}. In the next paper \cite{DS3} they showed that the value of the uniform exponent
for {\it simultaneous} approximations to any point $ (\theta,\theta^2)$ is not greater than $\frac{\sqrt{5}-1}{2}$. This together with  Jarn\'{\i}k's transference equality
(see \cite{JTb}) leads to the bound $\hat{\omega} \le \frac{3+\sqrt{5}}{2}$  which holds for all linear forms with coefficients of the form $ \theta,\theta^2$.
So for a linear form with coefficients   $ \theta,\theta^2$ one has
\begin{equation}\label{sec}
2\le \hat{\omega} \le \frac{3+\sqrt{5}}{2}
.
\end{equation}
D. Roy \cite{Ro1,Ro2} showed that the set of values $\hat{\omega}$ for linear forms under our consideration form a dense set in the segment (\ref{sec}).
Moreover  he constructed a countable set of numbers $\theta$ such that
$$
\hat{\omega}(\theta,\theta^2) = \frac{3+\sqrt{5}}{2}\,\,\,\,\,\text{and}\,\,\,\,\,
\omega_*(\theta) =  3+\sqrt{5}.
$$
This shows that our bound (\ref{t2}) from Theorem 2 is optimal in  the right endpoint of the segment (\ref{sec}), namely for $\hat{\omega} = \frac{3+\sqrt{5}}{2}$.

Other results on approximation by algebraic numbers are discussed in W. Schmidt's book \cite{SB},
in wonderful book by Y. Bugeaud \cite{BB} and in M. Waldschmidt's survey \cite{W}.

Our proof of Theorem 1 generalizes ideas from \cite{DS1,DS2,DS3} and uses Jarn\'{\i}k's inequalities \cite{JSZ,JRUS}.

{\bf 3. Minimal points.}

In the sequel we may suppose that $ \hat{\omega} >2$ as the case $\hat{\omega} = 2$ follows from Davenport-Schmidt's theorem (in this case
our Theorem 1 claims that $ \omega_{LP} \ge 3$).
We take $\alpha <\hat{\omega}$ close to $\hat{\omega}$ so that $ \alpha >2$.

A vector $ {\bf x } = (x_0,x_1,x_2) \in \mathbb{Z}^3\setminus \{{\bf 0}\}$ is defined to be a {\it minimal point} (or {\it best approximation}) if
$$
\min_{{\bf x}': \, 0< |{\bf x}'|\le |{\bf x}|} |L({\bf x}')|
= L({\bf x}).
$$
As $1,\theta_1,\theta_2$ are linearly independent, all
 the minimal points form a sequence ${\bf x}_\nu= (x_{0,\nu},x_{1,\nu},x_{2,\nu})$, $ \nu =1,2,3,...$ such that for $ X_\nu =|{\bf x}_\nu|, L_\nu = L({\bf x}_\nu)$ where one has
$$
X_1 < X_2<...<X_\nu< X_{\nu+1}<...\,\, ,\,\,\,\,\,
L_1>L_2>...> L_\nu > L_{\nu+1}>... \, .
$$
Here we should note that
\begin{equation}\label{sii}
L_j \le X_{j+1}^{-\alpha}
\end{equation}
 for all $j$ large enough.
Of course each vector ${\bf x}_j$ is  primitive and each  couple $ {\bf x}_j, {\bf x}_{j+1}$ form a basis of the two-dimensional lattice
$ \mathbb{Z}^3 \cap {\rm span}\, ({\bf x}_j, {\bf x}_{j+1})$.

Let $F({\bf x})$ be a linear form linearly independent with $ L$ and $P$. Then
\begin{equation}\label{max}
\max \{ |L({\bf x})|, |P({\bf x})| , |F({\bf x})|\}
\asymp |{\bf x}|.
\end{equation}
We also use the notation $ P_\nu = P({\bf x}_\nu), F_\nu = F({\bf x}_\nu)$.
In the sequel we need to consider determinants
$$
\Delta_j
=
\left|
\begin{array}{ccc}
L_{j-1}& P_{j-1}&F_{j-1}
\cr
L_{j}& P_{j}&F_{j}
\cr
L_{j+1}& P_{j+1}&F_{j+1}
\end{array}
\right|
=
A
\left|
\begin{array}{ccc}
x_{0.j-1}& x_{1,j-1}&x_{2,j-1}
\cr
x_{0,j}& x_{1,j}&x_{2,j}
\cr
x_{0,j+1}& x_{1,j+1}&x_{2,j+1}
\end{array}
\right|
,
$$
here $A$ is a non-zero constant depending on the coefficients of linear forms $L,P,F$.
We take into account (\ref{max},\ref{sii}) and the inequality $\alpha >2$ to see that
\begin{equation}\label{detta}
\Delta_j  = L_{j-1}P_{j} F_{j+1} - L_{j-1}P_{j+1}F_j +O(L_j X_{j+1}^2) =
L_{j-1}P_{j} F_{j+1} - L_{j-1}P_{j+1}F_j +o(1),\,\,\, j\to \infty.
\end{equation}
The following statement is a variant of Davenport-Schmidt's lemma.  We give it without a proof. It deals with three consecutive minimal points
${\bf x}_{j-1}, {\bf x}_j, {\bf x}_{j+1}$
 lying in a two-dimensional linear subspace, say $\pi$. We should note that our definition of minimal points differs from those in
\cite{DS1, DS2,SB}. However the main argument is the same. It is discussed in our survey \cite{MH}. One may look
for the  approximation of   the one dimensional subspace $ \ell =\pi \cap \{ {\bf z}:\, L({\bf z}) = 0\}$
by the points of two-dimensional lattice $\Lambda_j = \langle x_{j-1}, x_j\rangle$
Then  the points  ${\bf x}_{j-1}, {\bf x}_j, {\bf x}_{j+1}\in \Lambda_j$ are the consecutive best approximations to $\ell$ with respect to the {\it induced} norm
on $\pi$ (see \cite{MH}. Section 5.5).

{\bf Lemma 1.}\,\,{\it If
for some $j$ the points ${\bf x}_{j-1}, {\bf x}_j, {\bf x}_{j+1}$ are linearly dependent then
$$
{\bf x}_{j+1}= t {\bf x}_j+ {\bf x}_{j-1}
$$
for some integer $t$.}

The next statement is known for long time. It comes from Jarn\'{\i}k's papers \cite{JSZ,JRUS}. It was rediscovered by Davenport and Schmidt in \cite{DS3} and discussed in our survey \cite{MH}.

{\bf Lemma 2.}\,\,{\it
there exist infinitely many indices $j$ such that the vectors  ${\bf x}_{j-1}, {\bf x}_j, {\bf x}_{j+1}$
are linearly independent.
}

 The following lemma is due to  Jarn\'{\i}k \cite{JSZ, JRUS} (see also Section 5.3 from our paper \cite{MH}).

{\bf Lemma 3.}\,\,{\it
Suppose that $j$ is large enough and the points ${\bf x}_{j-1}, {\bf x}_j, {\bf x}_{j+1}$ are linearly independent. Then
\begin{equation}\label{ja}
X_{j+1} \gg X_j ^{\alpha - 1}
\end{equation}
and
\begin{equation}\label{ja1}
L_j\ll X_j ^{-\alpha(\alpha - 1)}
\end{equation}

}

Now we take large $\nu$ and $ k \ge \nu+1$ such that

$\bullet$ vectors ${\bf x}_{\nu-1}, {\bf x}_\nu, {\bf x}_{\nu+1}$
are linearly independent;

$\bullet$ vectors ${\bf x}_{k-1}, {\bf x}_k, {\bf x}_{k+1}$
are linearly independent;

$\bullet$ vectors $ {\bf x}_j, \, \nu \le j \le k$ belong to the two-dimensional lattice
$\Lambda_\nu = \mathbb{Z}^3 \cap {\rm span} \,({\bf x}_\nu, {\bf x}_{\nu+1})$.

From Lemma 1  it follows that
for $j$ from the range $\nu\le j \le k-1$ one has
$$
L_{j+1} = t_{j+1} L_j + L_{j-1},\,\,\,\,\, P_{j+1} = t_{j+1} P_j + P_{j-1},
$$
with some integers $t_{j+1}$, and hence
\begin{equation}
\label{iii}
L_\nu P_{\nu+1}- L_{\nu+1}P_\nu = \pm (L_{k-1}P_k+L_kP_{k-1}).
\end{equation}

{\bf Lemma 4.}\,\, {\it Consider  positive $r$ under the condition
\begin{equation}\label{er}
r < \alpha^2-\alpha+1 < \hat{\omega}^2 - \hat{\omega}+1.
\end{equation}
Suppose that
\begin{equation}\label{r1}
|P_\nu| \le L_\nu X_\nu^r
\end{equation}
and $\nu$ is large. Then
\begin{equation}\label{llee}
|P_{\nu+1}| \gg X_\nu^{\alpha - 1}.
\end{equation}}

Proof.
For $ j = \nu$
consider the second term in the r.h.s of (\ref{detta}). From (\ref{sii},\ref{max},\ref{er}.\ref{r1}) and the inequality (\ref{ja}) of Lemma 3  we have
$$
|L_{\nu-1}P_\nu F_{\nu+1}| \ll |L_{\nu-1} L_\nu X_\nu^r|\,X_{\nu+1}\ll
X_\nu^{r-\alpha}X_{\nu+1}^{1-\alpha} \ll X_{\nu}^{r -\alpha^2+\alpha - 1} = o(1).
$$
As $\Delta_\nu \neq 0$ we see that
$$
1\ll |L_{\nu-1} P_{\nu+1} F_\nu| \ll L_{\nu-1} | P_{\nu+1}|\, X_\nu \ll X_{\nu}^{1-\alpha}| P_{\nu+1}|
$$
(in the last inequalities we use (\ref{max}) and (\ref{sii}). Everything is proved.$\Box$.

{\bf 4. The main estimate.}

The following Lemma  presents our main argument.

{\bf Lemma 5.}\,\,{\it
Suppose that $r$ satisfies (\ref{er}).
 Suppose that (\ref{sii}) holds for all indices $j$ and  suppose that for a certain
 $\beta_0$ one has
\begin{equation}\label{or}
L_\nu\gg X_\nu^{-\beta_0}
.
\end{equation}
 Suppose that simultaneously we have
\begin{equation}\label{a1}
|P_{\nu}| \le L_{\nu} X_{\nu}^r,
\end{equation}
\begin{equation}\label{a2}
|P_{k-1} |\le L_{k-1} X_{k-1}^r,
\end{equation}
\begin{equation}\label{a3}
|P_{k}| \le L_{k} X_\nu^r.
\end{equation}
Then
\begin{equation}\label{itog1}
r \ge \alpha^2 + 1 - \frac{\beta_0}{\alpha - 1}.
\end{equation}
and
\begin{equation}\label{itog2}
L_k \gg X_k^{-\beta'},\,\,\,\,\text{with}\,\,\,\, \beta' = r-\alpha-1+\frac{\beta_0}{\alpha-1}<\beta_0.
\end{equation}
}

 First of all we note that
$$
L_{\nu+1}|P_\nu|\le
L_\nu L_{\nu+1}X_\nu^r \ll
L_\nu X_{\nu+2}^{-\alpha} X_\nu^r \ll
L_\nu X_{\nu+1}^{-\alpha} X_\nu^r\ll
L_\nu X_\nu^{r-\alpha(\alpha-1)}= o( L_\nu X_\nu^{\alpha -1}).
$$
Here the first inequality comes from (\ref{a1}). The second inequality is (\ref{sii}  with $j = \nu+1$.
The third one is simply $ X_{\nu+2}\ge X_{\nu+1}$.
The fourth one is (\ref{ja}) of Lemma 3 for $ j = \nu$.
The last inequality here follows from (\ref{er}) as
$ r <\alpha^2- \alpha + 1 <\alpha^2 -1$ (because $\alpha >2$).
We see that the conditions of Lemma 4 are satisfies and by Lemma  4 we see that
$$ L_\nu |P_{\nu+1} | \gg L_\nu X_\nu^{\alpha - 1}.$$
So in the l.h.s. of (\ref{iii}) the first summand is larger than the second.
Now from (\ref{iii}) we have
\begin{equation}\label{au1}
L_\nu X_\nu^{\alpha - 1}\ll L_{k-1}|P_k|+L_k|P_{k-1}|.
\end{equation}
We apply (\ref{a2},\ref{a3}) to see that
\begin{equation}\label{au2}
\max ( L_{k-1}|P_k|,L_k|P_{k-1}|) \le L_{k-1}L_k X_k^r\ll X_k^{r-\alpha} X_{k+1}^{-\alpha} \le X_k^{r-\alpha^2}\le
X_{\nu+1}^{r-\alpha^2}\ll X_\nu^{(r-\alpha^2)(\alpha-1)}.
\end{equation}
Here the second inequality comes from (\ref{iii}) for $ j = k-1$ and $j = k$.
The third inequality is Lemma 3 with $j = k$. The fourth one is just $ X_k \ge X_{\nu+1}$.
The fifth one is Lemma 3 for $j=\nu$.

Now from estimates (\ref{au1},\ref{au2}) and (\ref{or}) we have
$$
X_\nu^{-\beta_0+\alpha - 1}\ll X_\nu^{(r-\alpha^2 )(\alpha - 1)}.
$$
This gives
$$
r \ge \alpha^2 + 1 - \frac{\beta_0}{\alpha - 1}.
$$
So (\ref{itog1}) is proved.

To get (\ref{itog2}) we combine the estimate (\ref{au1}) with the left inequality of (\ref{au2}), the bound (\ref{or}) for $ j = \nu$ and  the bound (\ref{sii}) for $ j = k-1$.
This gives
$$
X_{\nu}^{\alpha-1-\beta_0}\le
L_\nu X_\nu^{\alpha - 1}\ll
L_{k-1}L_k X_k^r\ll L_k X_k^{r-\alpha},
$$
or
$$
L_k \gg  X_k^{\alpha -r} X_\nu^{\alpha-1-\beta_0}.
$$
But $\beta_0 > \alpha(\alpha-1) \ge (\alpha-1)$ by inequality (\ref{ja1}) of Lemma 3 and $X_k \ge X_{\nu+1}\gg X_\nu^{\alpha-1}$ by
inequality (\ref{ja}) of Lemma 3. So
$$
L_k \gg X_k^{\alpha-r +\frac{\alpha -1-\beta_0}{\alpha -1}},
$$
and this is  the first inequality form (\ref{itog2}).

Moreover as $\beta_0 >\alpha (\alpha -1) $, from (\ref{er}) we deduce
$ \beta'<\beta$. Lemma is proved.$\Box$

{\bf 5. Proof of Theorem 1.}

Suppose that $r$ satisfies (\ref{er}).
We take infinite sequence  indices $ \nu_1<\nu_2<...<\nu_t<...$ such that

$\bullet$
for every $i= 1,2,... $
vectors ${\bf x}_{\nu_i-1}, {\bf x}_{\nu_i}, {\bf x}_{\nu_i+1}$
are linearly independent;

$\bullet$ for $i = 1,2,... $ vectors $ {\bf x}_j, \, \nu_i \le j \le \nu_{i+1}$ belong to the two-dimensional lattice
$\Lambda_{\nu_i} = \mathbb{Z}^3 \cap {\rm span} \,({\bf x}_{\nu_i}, {\bf x}_{\nu_i+1})$.

Now we suppose that three inequalities (\ref{a1},\ref{a2},\ref{a3}) hold for all  triples
$(\nu,k-1,k) = (\nu_i,\nu_{i+1}-1,\nu_{i+1})$ for all $i \ge 1$.

Define recursively
$$
\beta_{i+1} = r-\alpha-1+\frac{\beta_i}{\alpha-1}.
$$
Then
$$
\beta_i = \alpha (\alpha -1) +\frac{\beta_0}{(\alpha-1)^i} \to \alpha (\alpha -1),\,\,\,\,\, i\to \infty .
$$
We apply of Lemma 5 to the first $w$ triple of indices. Then we get  (\ref{itog2}) for $k = \nu_{i+1}$, and in particular for $ k=\nu_w$
with
$
\beta_w
$
close to $\alpha (\alpha - 1)$.
Now we apply Lemma 5 to $\nu = \nu_w$. in (\ref{or}) we have $\beta_w$ instead of $\beta$.
So (\ref{itog1}) gives
$$
r \ge \alpha^2 - \alpha + 1 - \frac{\beta_w}{\alpha - 1}
 $$
We take limit $ w\to \infty$ to see that
$$
r \ge \alpha^2 - \alpha+1.
$$

This contradicts to (\ref{er}).
So there exists $j \in \cup_{i=1}^\infty \{ \nu_i, \nu_{i+1}-1,\nu_{i+1}\}$ such that
$L_j\le |P_j|X_j^{-r}$. Theorem is proved.$\Box$

{\bf Acknowledgement.} The author is grateful to Igor Rochev
for important comments on the proof and for pointing out certain inaccuracies in the manuscript.

\end{document}